\documentclass[11pt,a4paper]{article}
\usepackage[margin=28mm]{geometry}
\usepackage[T1]{fontenc}

\usepackage{lmodern,amsmath,amssymb,amsthm,microtype,enumitem}
\usepackage[colorlinks=true,linkcolor=blue,citecolor=blue,urlcolor=blue]{hyperref}
\usepackage{comment}
\usepackage{xcolor}
\usepackage{stmaryrd,mathrsfs,bm,bbm,mathtools}

\newtheorem{theorem}{Theorem}[section]
\newtheorem{proposition}[theorem]{Proposition}
\newtheorem{lemma}[theorem]{Lemma}
\newtheorem{corollary}[theorem]{Corollary}
\newtheorem{definition}[theorem]{Definition}
\newtheorem{remark}[theorem]{Remark}
\numberwithin{equation}{section}

\newcommand{\HH}{\mathcal H}
\newcommand{\LL}{\mathcal L}
\newcommand{\R}{\mathbb R}
\newcommand{\D}{\mathbb D}

\renewcommand{\S}{\mathbb S}
\newcommand{\md}{\operatorname{md}}
\newcommand{\apmd}{\operatorname{ap\,md}}
\newcommand{\Mod}{\operatorname{Mod}}
\newcommand{\loc}{\mathrm{loc}}
\newcommand{\J}{\mathbf{J}}
\newcommand{\ud}{\operatorname{d}}
\newcommand{\aplim}{\operatorname{ap}\lim}

\title{Infinitesimally Euclidean metric surfaces}

\author{Alexander Lytchak, Elefterios Soultanis}

\date{\today}
\begin{document}
\maketitle

\begin{abstract}
We prove the equivalence of several notions of infinitesimally Euclidean-like behaviour for metric surfaces appearing in the literature. We obtain in particular that metric surfaces with property (ET) are universally infinitesimally Hilbertian.
\end{abstract}

\section{Introduction}
In this short note, we prove the equivalence of three notions describing ``Euclidean behaviour at infinitesimal scales'' for metric surfaces, i.e. metric spaces $(X,d)$ homeomorphic to a surface without boundary for which the Hausdorff 2-measure $\mathcal H^2$ is locally finite. The notions we consider are infinitesimal Hilbertianity \cite{Gigli15}, property \textup{(ET)} \cite[Section~11]{LW17}, and having Euclidean approximate tangents \cite[Section~5]{AK}.

Infinitesimal Hilbertianity was introduced by Gigli \cite{Gigli15} and is a key part of the definition of RCD-spaces \cite{AGS14}; it has also attracted independent interest \cite{Ryborz26,Ryborz26b,ES,ERS24,ERS24a,NZPS,APS15,APS19,DMGPS21,LP20}. Property (ET), defined by Lytchak and Wenger \cite[Section~11]{LW17}, arose in the study of isoperimetric inequalities on metric spaces. In particular, the Euclidean isoperimetric inequality for curves implies property (ET), and this implication is used in isoperimetric characterizations of upper curvature bounds; see \cite{LW18} and \cite[Section~5]{SW25}. Approximate tangents are closely related to rectifiability of metric spaces; see \cite{K94,AK,Bate} for related results.

A metric length surface is a metric surface that is also a length space.

\begin{theorem}\label{thm:main}
Let $X$ be a metric length surface. Then the following are equivalent:
\begin{enumerate}[label=\textup{(\roman*)}]
\item $X$ has property \textup{(ET)}.
\item $(X,d,\mathcal H^2)$ is infinitesimally Hilbertian.
\item Every countably $2$-rectifiable subset of $X$ has a Euclidean
approximate tangent plane in the sense of Ambrosio--Kirchheim at $\mathcal H^2$-almost every point.
\end{enumerate}
\end{theorem}
For the definition of infinitesimal Hilbertianity, property (ET), and approximate tangent planes we refer to Section~2.

Theorem \ref{thm:main} justifies \cite[Remark~1.2]{LS26}, and its proof is based on the uniformization results of Ntalampekos--Romney and Meier--Wenger \cite{NR26,MW25}, together with the analysis in \cite{MR} of the \emph{singular part} of the surface, that is the image of the region where a weakly quasiconformal parametrization fails the Lusin property (N); see Definition \ref{def:sing-part}.

The length assumption is needed only for the implication (i)$\implies$(iii): the converse (iii)$\implies$(i) holds for arbitrary metric spaces; see Subsection~\ref{sec:general-relations}. Metric length surfaces need not be rectifiable; instead we prove in Proposition \ref{prop:unrect} that their singular part is purely 2-unrectifiable. This need not be the case without the length assumption \cite[Theorem~5.3]{EIR22}.

For general metric surfaces $(X,d)$ equipped with $\HH^2$, the singular part nevertheless forms a zero-dimensional chart in the sense of \cite{ES}. In Lemma \ref{lem:no-1-dim-charts} we prove that in metric surfaces every zero-dimensional chart is contained in the singular part up to a $\HH^2$-null set, and all other charts are 2-dimensional. This leads to the following theorem.

\begin{theorem}\label{thm:main2}
If $X$ is a metric surface, then the following are equivalent:
\begin{enumerate}[label=\textup{(\roman*)}]
\item $X$ has property (ET);
\item $(X,d,\mathcal{H}^2)$ is infinitesimally Hilbertian.
\end{enumerate}
\end{theorem}

The measure $\HH^2$ plays no special role when it comes to infinitesimal Hilbertianity; in fact, either of the two conditions in Theorem \ref{thm:main2} is equivalent to the metric surface $X$ being universally infinitesimally Hilbertian: 

\begin{theorem}\label{thm:main3}
Suppose $X$ is a metric surface satisfying (i) or (ii) in Theorem \ref{thm:main2}. Then $(X,d,\mu)$ is infinitesimally Hilbertian for any locally finite Borel regular measure $\mu$. 
\end{theorem}
Theorem \ref{thm:main3} gives new examples of universally infinitesimally Hilbertian spaces; previous (non-smooth) examples rely on curvature bounds \cite{NZPS,DMGPS21} or a sub-Riemannian structure \cite{LP20,LDLP23}, whereas Theorem \ref{thm:main3} relies on 2-dimensionality. In particular, by a recent result of Ryborz \cite{Ryborz26b}, 2-manifolds $(M,d_g)$ with (possibly discontinuous) Riemannian metric $g$ with $g,g^{-1}\in L^\infty_\loc\cap W^{1,p}_\loc$, $p\geq1$, are universally infinitesimally Hilbertian.

\subsubsection*{Acknowledgements}
The authors would like to thank Vanessa Ryborz and Timo Schultz for helpful exchanges on closely related subjects. The first author was supported by the DFG under Germany’s Excellence Strategy -- EXC-2047/1 -- 390685813. The second author was supported by the Research Council of Finland grant no. 355122.
 
This paper arose in close collaboration with ChatGPT. The collaboration involved exploring and refining arguments, locating references, and drafting and revising the manuscript. The authors take full responsibility for the mathematical content and the final text.

\section{Preliminaries}
\subsection{Definitions and notation}
Let $X$ and $Y$ be separable metric spaces, $\mu$ a locally finite Borel regular measure on $X$, and $u:X\to Y$ a mapping. A Borel function $g:X\to [0,\infty]$ is called an \emph{upper gradient} of $u$, if 
\begin{align}\label{eq:ug}
	d_Y(u(\gamma(b)),u(\gamma(a)))\le \int_\gamma g\ud s
\end{align}
for all absolutely continuous curves $\gamma:[a,b]\to X$. If \eqref{eq:ug} holds for all $\gamma\notin \Gamma_0$, where $\Gamma_0$ is a family of curves satisfying
\begin{align*}
	\Mod_2(\Gamma_0):=\inf\Big\{\int\rho^2\ud\mu:\  \int_\gamma\rho\ud s\ge 1\ \forall \gamma\in \Gamma_0  \Big\}=0,
\end{align*}
we say that $g$ is a weak upper gradient of $u$. We refer to Bj\"orn--Bj\"orn \cite{BB11} and Heinonen--Koskela--Shanmugalingam--Tyson \cite{HKST15} for the definition of absolutely continuous curves and the curve integral, as well as further discussion. The (local) Newton--Sobolev space $W^{1,2}_{(\loc)}(X,Y)$ consists of maps $u:X\to Y$ with $u\in L^2_{(\loc)}(X,Y)$ for which there exists a locally 2-integrable upper gradient $g:X\to [0,\infty]$.\footnote{In the literature, the notation $N^{1,2}_\loc(X,Y)$ is more common.} If $u\in W^{1,2}_{(\loc)}(X,Y)$, there exists a pointwise minimal weak upper gradient in $L_{(\loc)}^2(\mu)$, denoted by $|Du|$. 

A metric measure space $(X,d,\mu)$ is called \emph{infinitesimally Hilbertian} if the space $W^{1,2}(X):=W^{1,2}(X,\R)$ is a Hilbert space with the norm $(\|u\|_{L^2(\mu)}^2+\||Du|\|_{L^2(\mu)}^2)^{1/2}$. A metric space $X$ is called \emph{universally infinitesimally Hilbertian} if $(X,d,\mu)$ is infinitesimally Hilbertian for every locally finite Borel regular measure $\mu$ on $X$.

A subset $E$ of a metric space $X$ is called \emph{countably $n$-rectifiable} if
\[
E=N\cup\bigcup_{i=1}^{\infty}u_i(K_i),
\]
where $\HH^n(N)=0$, the sets $K_i\subset\R^n$ are compact, and the
maps $u_i:K_i\to X$ are bi-Lipschitz. For $\LL^n$-almost every
$z\in K_i$ there is a metric differential, that is, a norm
$\md_z u_i$ on $\R^n$ satisfying
\[
\lim_{\substack{x,y\to z\\x,y\in K_i}}
\frac{d_X(u_i(x),u_i(y))-\md_z u_i(x-y)}{|x-z|+|y-z|}=0;
\]
see \cite{K94}. For $\HH^n$-almost every $x\in E$, the isometry class
of $(\R^n,\md_z u_i)$, where $x=u_i(z)$, is independent of the
chosen rectifiable decomposition and parametrization \cite[Section~5]{AK}.
It is called the \emph{approximate tangent space} of $E$ at $x$.
We say that $E$ has Euclidean approximate tangents if these norms are
induced by inner products at $\HH^n$-almost every point of $E$.

Sobolev spaces on Euclidean domains, in particular on the unit disc $\D$, are always taken with respect to Lebesgue measure. For any $u\in W^{1,2}(\D,Y)$, there exists a seminorm $\apmd_zu$ on $\R^2$ for a.e. $z\in \D$, satisfying 
\begin{align}
	\aplim_{x,y\to z}\frac{d_Y(u(x),u(y))-\apmd_zu(x-y)}{|y-z|+|x-z|}=0.
\end{align}
In this case, the minimal weak upper gradient is given by $|Du|(z)=\sup\{\apmd_zu(v): v\in \S^1\}$ $\mathcal L^2$-a.e. $z\in\D$; see \cite[Proposition~4.3 and Section~2]{LW17} and \cite[Theorem~7.1.20]{HKST15}.   The space $X$ has \emph{property (ET)} if for every $u\in W^{1,2}(\mathbb D,X)$ the seminorm $\apmd _z u$ is the square root of a positive semidefinite quadratic form at almost every $z\in \mathbb D$.

For a seminorm $s$ on $\R^2$, define its (Hausdorff) \emph{Jacobian} by
\[
\J(s)=\frac{\pi}{\LL^2(\{v\in\R^2:s(v)\leq1\})},
\]
with value zero when the denominator is infinite. Thus $\J(s)=0$ exactly when $s$ is degenerate. We write $\J f(z):=\J(\apmd_zf)$ wherever the approximate metric differential exists.

When $Y=\R$ we denote $W^{1,2}_\loc(X,\R)=W^{1,2}_\loc(X)$. In \cite{ES} it was proved that if $\dim_HX<\infty$ and $X$ is locally compact, then $X$ can be covered (up to a $\mu$-null set) by a countable collection $\{(U_i,\varphi_i)\}$ of 2-weak charts, that is $U_i\subset X$ is Borel with $\mu(U_i)>0$ and $\varphi_i:X\to \R^{n_i}$ is a Lipschitz map for each $i$, characterized by the following property: for any Lipschitz function $f:X\to \R$, there exists a pointwise $\mu$-a.e. unique Borel function $\ud f:U_i\to (\R^{n_i})^*$ 
for which
\begin{align}\label{eq:diff}
	(f\circ\gamma)'(t)=\ud_{\gamma(t)}f((\varphi_i\circ\gamma)'(t))\quad\textrm{a.e.}\ t\in \gamma^{-1}(U_i)
\end{align}
for all absolutely continuous curves $\gamma$ outside some curve family of zero 2-modulus. In this paper, the term \emph{chart} means a 2-weak chart, and the dimension of a chart is the dimension of the target space of the chart mapping. If $(U,\varphi)$ is an $n$-dimensional chart, then $n\le \dim_H U$ \cite[Proposition~4.13]{ES}. Note that $(U,\varphi)$  can be an $m$-dimensional chart with $m<n$ even if $U$  is bi-Lipschitz to a subset of $\R^n$  and $\mathcal H^n (U) >0$.

When $\varphi\in W^{1,2}_\loc(X,\R^N)$, there exists a \emph{canonical minimal gradient}, that is a Borel map $\Phi:X\times(\R^N)^*\to [0,\infty]$ characterized by the following:
\begin{enumerate}[label=\textup{(\roman*)}]
	\item For each $\xi\in (\R^N)^*$, $\Phi_\xi:=\Phi(\cdot,\xi)$ is a representative of $|D(\xi\circ \varphi)|$;
	\item for $\mu$-a.e. $x\in X$, $\Phi^x:=\Phi(x,\cdot)$ is a seminorm on $(\R^N)^*$. 
\end{enumerate}

A pair $(U,\varphi)$, where $U\subset X$ is Borel with $\mu(U)>0$ and $\varphi:X\to\R^n$ is Lipschitz, is a chart of dimension $n$ if and only if $\Phi^x$ is a norm on $(\R^n)^*$ for $\mu$-a.e. $x\in U$ and the dimension is maximal, i.e. if $\psi:X\to \R^m$ is Lipschitz and $m>n$, then the canonical minimal gradient $\Psi$ associated with $\psi$ satisfies that $\Psi^x$ is degenerate $\mu$-a.e. $x\in U$. We denote the dual norm $(\Phi^x)^*=|\cdot|_x$ and call it the pointwise norm of the chart. It is independent of the chart in the sense that the pointwise norms of any two charts overlapping in a set of positive $\mu$-measure are linearly isometric for $\mu$-a.e. point of their common domain of definition. Moreover, $|Df|(x)=|\ud_x f|_x^*$ for $\mu$-a.e. $x\in U$ for any $f\in W^{1,2}_\loc(X)$. See \cite{ES,ERS24} for these facts. In dimension zero the pointwise tangent space is trivial; a Borel set carries a zero-dimensional chart precisely when the minimal weak upper gradients of all Lipschitz functions vanish there almost everywhere, cf. \cite[Proposition~4.4]{ES}.

\begin{remark} \label{rem:rem}
	In \cite{ES,ERS24} it is assumed that $X$ is complete. Since minimal upper gradients and the pointwise norms are local, the results there readily generalize to locally complete metric measure spaces, including metric surfaces.
\end{remark}

\subsection{Relations in general metric spaces}\label{sec:general-relations}
For a metric measure space admitting a finite-dimensional $2$-weak
differentiable structure, \cite[Corollary~6.7(ii)]{ES} gives
\begin{equation}\label{eq:ih-charts}
	(X,d,\mu)\text{ is infinitesimally Hilbertian}
	\Longleftrightarrow
	|\cdot|_x\text{ is an inner product norm for $\mu$-a.e. }x.
\end{equation}
Equivalently, the pointwise parallelogram identity holds for all  $f,g\in W^{1,2} (X)$: 
\[
|D(f+g)|^2+|D(f-g)|^2=2|Df|^2+2|Dg|^2
\quad\mu\text{-a.e.}
\]

In general neither  property (ET) nor infinitesimal Hilbertianity imply each other:  The Heisenberg group $(\mathbb H^1,d_{|\cdot|},\mathcal H^4 )$ with a non-Riemannian sub-Finsler  metric is not infinitesimally Hilbertian \cite[Section~1.5]{MagnaboscoRossi}, but has property (ET), since $\mathbb H^1$ does not contain $2$-rectifiable subsets of positive $\mathcal H^2$-measure \cite[Theorem~7.2]{AK}.  On the other hand, a non-Riemannian Finsler plane
does not have property (ET), but is infinitesimally Hilbertian with respect to a Dirac measure \cite[Introduction]{DMGPS21}.

The following result for  complete metric  spaces  has been obtained in \cite[Theorem~A.1]{NZPS}:  
\begin{equation}
	X\text{ is universally infinitesimally Hilbertian}  \Longrightarrow
	X \text{ has property (ET)}.
\end{equation}

On the other hand, for an arbitrary metric
space $X$, if every countably $2$-rectifiable subset of $X$ has Euclidean
approximate tangents almost everywhere, then $X$ has property  (ET). Indeed,
for $u\in W^{1,2}(\D,X)$, the set where $\apmd_z u$ is nondegenerate
can be covered, up to a null set, by compact bi-Lipschitz pieces
\cite[Proposition~4.3]{LW17}; on these pieces the approximate
tangent norm agrees with $\apmd_z u$ almost everywhere. Degenerate seminorms on
$\R^2$ are always induced by positive semidefinite quadratic forms.

\section{The proofs}

\subsection{A technical lemma}
The key technical lemma of this note expresses the pointwise norm in terms of the metric differentials of a Sobolev parametrization. The general form of the lemma stated below, with a measure $\mu\ge \HH^2$, will be needed for universal infinitesimal Hilbertianity; see Section \ref{sec:universal}.
	
\begin{lemma}\label{lem:key-lemma}
Let $X$ be a locally compact metric space so that $\HH^2$ is locally finite. Let $\mu$ be a locally finite Borel regular measure on $X$ satisfying $\HH^2_X\le\mu$. Suppose $u\in W^{1,2}(\D,X)$ and $K\subset\D$ is a Borel set such that $u|_K$ is bi-Lipschitz. If $\varphi:X\to \R^2$ is a Lipschitz extension of $u^{-1}|_{u(K)}$, then for $\HH^2$-a.e. $x\in u(K)$
\begin{align}\label{eq:key-lemma}
(\Phi^x)^*(v)=\apmd_{\varphi(x)}u(v) ,\quad v\in\R^2.
\end{align}
Here $\Phi$ is the canonical minimal gradient of $\varphi$ with respect to $\mu$.
\end{lemma}

\begin{remark}
In Lemma \ref{lem:key-lemma}, we do not claim that \eqref{eq:key-lemma} holds for $\mu$-a.e. $x\in u(K)$. Indeed, this need not be true as $(X,d,\mu)$ can have e.g. 1-dimensional charts. The assumption $\HH^2\le\mu$ cannot be omitted in the statement for similar reasons.
\end{remark}

If $u$ is weakly quasiconformal and has essentially bounded multiplicity, the lemma also follows by adapting the argument in \cite[Appendix A]{NZPS}. The proof below requires neither assumption. Lemma \ref{lem:key-lemma} has the following corollary, which will be important  in the sequel.
\begin{corollary}\label{cor:rect-charts}
	Under the hypotheses of Lemma \ref{lem:key-lemma}, if $\LL^2(K)>0$, we have that $(u(K), \varphi)$ is a 2-dimensional chart of the metric measure space $(X,d,\HH^2)$. 
\end{corollary}
\begin{proof}
Apply Lemma \ref{lem:key-lemma} with $\mu=\HH^2$. By \eqref{eq:key-lemma} and \cite[Lemma 4.3(1)]{ES}, $(u(K),\varphi)$ is 2-independent, while by \cite[Proposition 4.13]{ES} the dimension of $(u(K),\varphi)$ is maximal. Thus $(u(K),\varphi)$ is a 2-dimensional chart.
\end{proof}

{
\begin{proof}[Proof of Lemma \ref{lem:key-lemma}]
Write $E=u(K)$ and $s_z=\apmd_zu$. We fix a seminorm representative  $\Phi^x\le\operatorname{Lip}(\varphi)|\cdot|$  everywhere. Set $|\cdot |_x = (\Phi^x)^{\ast}$  allowing the dual to take the value $\infty$. Note that all weak upper gradients and moduli in this proof are taken with respect to $\mu$. 

We can exhaust $K$, up to a $\LL^2$-null set, by compact subsets. Since $u|_K$ is Lipschitz, their images exhaust $E$ up to a $\HH^2$-null set. It thus suffices to assume that $K$ is compact. 
 
We find a Borel subset $N \subset E$ with $\mathcal H^2(N)=0$, and such that for all $z\in K\setminus \varphi (N)$ the metric differential of $u|_K$ at $z$ coincides with $s_z$, \cite[Proposition 4.3]{LW17}.   Hence, $s_z$ is a norm. For every absolutely continuous curve $\gamma$ and a.e. $t\in\gamma^{-1}(E\setminus N)$, metric differentiability and differentiation along this measurable set give:
\begin{equation}\label{eq:curvewise-metric}
|\gamma'|(t)=s_{\varphi(\gamma(t))}((\varphi \circ \gamma ) '(t))\,.
\end{equation}
There is a curve family $\Gamma_0$ with $\Mod_{2,\mu}(\Gamma_0)=0$ such that all $\gamma\notin\Gamma_0$ satisfy
\begin{equation} \label{eq:curvewise-tangent}
|(\varphi \circ \gamma )' (t)|_{\gamma(t)}\le |\gamma'|(t)
\quad\text{for a.e. }t.
\end{equation}
To construct $\Gamma_0$, take a countable dense family of $\xi\in(\R^2)^*$ and include, for each $\xi$, all curves
violating  the scalar weak upper gradient inequalities of $\xi\circ\varphi$ on some subcurve.   For every remaining $\gamma$, differentiation gives, simultaneously for these $\xi$,
\[
|\xi((\varphi \circ \gamma)'(t))|\le\Phi^{\gamma(t)}(\xi)|\gamma'|(t)
\quad\text{for a.e. }t.
\]
Continuity in $\xi$ extends this to all $\xi\in(\R^2)^*$, and taking duals proves \eqref{eq:curvewise-tangent}.

Fix $\xi\in(\R^2)^*$ and set $h=\xi\circ\varphi$. 
For every $\gamma\notin\Gamma_0$, \eqref{eq:curvewise-metric} and duality yield
\[
|(h\circ\gamma)'(t)|=|\xi((\varphi \circ \gamma)'(t))|\le s_{\varphi(\gamma(t))}^*(\xi)|\gamma'|(t)
\quad\text{a.e. on }\gamma^{-1}(E\setminus N).
\]
On the complement of $E\setminus N$ in $X$, \eqref{eq:curvewise-tangent} gives
$|(h\circ\gamma)'(t)|\le\Phi^{\gamma(t)}(\xi)|\gamma'|(t)$.
Thus the function
\[
 s_{\varphi(\cdot)}^*(\xi)\chi_{E\setminus N}
 +\Phi^{\cdot}(\xi)\chi_{X\setminus(E\setminus N)}
\]
is a  weak upper gradient of $h$. Note that the function is bounded above, hence locally square-integrable since $\mu$ is locally finite. By local minimality, this function is at least $|Dh|$ $\mu$-almost everywhere, and hence $\HH^2$-almost everywhere, on $X$.  Hence, for $\LL^2$-a.e. $z\in K$,  
$|Dh|(u(z)) \leq s^{\ast}_z (\xi)\,.$
 Taking a countable dense set of $\xi$ and dualizing gives $s_z\le |\cdot|_{u(z)}$ for $\LL^2$-a.e. $z\in K$.

For the reverse inequality, fix $v\in\S^1$, choose a unit vector
$e\perp v$, and put $\psi=e\cdot\varphi$, $F=\psi\circ u$.    
Choose an \(L^2\) upper gradient \(G\) of \(u\), a relatively compact open neighbourhood \(V\) of \(E\), and a nonnegative Borel function $\rho_0\in L^2(X,\mu)$ whose integral is infinite along every curve in $\Gamma_0$. Define $\rho=\rho_0$ on $X\setminus N$ and $\rho=\infty$ on $N$. Although $\rho$ need not belong to $L^2(\mu)$, we have
\[
\int_V\rho^2\,\ud\HH^2=\int_V\rho_0^2\,\ud\HH^2
\le\int_X\rho_0^2\,\ud\mu<\infty.
\]
In particular, $\rho\in L^1(V,\HH^2)$, and its integral is still infinite along every curve in $\Gamma_0$.

The Sobolev coarea formula  is applicable to the function $F$ (which is $2$-quasicontinuous, cf.  \cite[Theorem 7.4.5]{HKST15}), as explained in \cite[p.481]{MSZ}.  Therefore,
 \cite[Theorem~1.1]{MSZ} and the Lipschitz coarea inequality
\cite[Theorem~2.1]{MR} give
\begin{equation}\label{eq:two-coarea}
\int_\R\int_{F^{-1}(a)}G\,\ud\HH^1\,\ud a
+\int_\R\int_{V\cap\psi^{-1}(a)}\rho\,\ud\HH^1\,\ud a<\infty.
\end{equation}
Fix a closed disk $B$ with $K\subset\operatorname{int}B$ and $B\Subset\D$. The planar curve-decomposition theorem \cite[Corollary 1]{ACMM01} together with \cite[p.173]{Fusco} implies that, for almost every $a\in \R$, the fiber $F^{-1}(a)\cap \operatorname{int}B$ agrees, up to an \(\mathcal H^1\)-null set, with a countable union of injective rectifiable arcs. We restrict to such values $a$ for which, in addition,  both inner integrals in \eqref{eq:two-coarea} are finite.

Parametrize each of these arcs $\sigma$ by arclength.
By our choice of $a$, we have $\int_\sigma G\,\ud s<\infty$.
Since $G$ is an upper gradient of $u$, the curve
$\gamma=u\circ\sigma$ is absolutely continuous.
Moreover, $\psi\circ\gamma=F\circ\sigma=a$ almost everywhere,
and hence $\psi\circ\gamma\equiv a$ by continuity.

Every injective rectifiable arc in $V\cap\psi^{-1}(a)$ has
finite $\rho$-integral and hence lies outside $\Gamma_0$.
Since upper gradient inequalities need only be verified on
injective arcs, 
\eqref{eq:curvewise-tangent} holds true along every absolutely continuous
curve in  $V\cap\psi^{-1}(a)$.

Since $\rho=\infty$ on $N$, our choice of $a$ also gives
$\HH^1(N\cap\psi^{-1}(a))=0$. As $\varphi$ is Lipschitz,
\[
\HH^1(\varphi(N)\cap F^{-1}(a))=0.
\]
For almost every $t\in\sigma^{-1}(K\setminus\varphi(N))$,
differentiating the identities
\[
\varphi\circ\gamma=\sigma,
\qquad e\cdot\sigma=a,
\]
which hold on this set, gives
\[
(\varphi\circ\gamma)'(t)=\sigma'(t),
\qquad e\cdot\sigma'(t)=0.
\]
Since $\sigma$ is parametrized by arclength, we have
$\sigma'(t)=\pm v$. Thus, writing $z=\sigma(t)$, metric
differentiability of $u|_K$ and \eqref{eq:curvewise-tangent} yield
\[
|v|_{u(z)}\le |\gamma'|(t)=s_z(v).
\]

The arcs cover $K\cap F^{-1}(a)$ up to an $\HH^1$-null set.
Hence this inequality holds for $\HH^1$-almost every point of
$K\cap F^{-1}(a)$, for almost every $a$.
Since $F(z)=e\cdot z$ on $K$, Fubini's theorem implies that it
holds for $\LL^2$-almost every $z\in K$.
Applying this to a countable dense set of $v\in\S^1$ and using
lower semicontinuity of $|\cdot|_{u(z)}$ gives
$|\cdot|_{u(z)}\le s_z$.
Together with the opposite inequality verified above, this proves
\eqref{eq:key-lemma}.
\end{proof}
}

\subsection{Metric surfaces}\label{sec:metr-surf}
{In this subsection we assume that $X$ is a metric surface. We use a countable cover by neighbourhoods $V\subset X$ admitting local Sobolev uniformizations $u:\D\to V$ of multiplicity one $\HH^2$-a.e. \cite[Theorem~2.8]{MR}.
	
\begin{definition}\label{def:sing-part}
	The \emph{singular part} $X_0$ of $X$ is characterized, up to $\HH^2$-null sets, by $X_0\cap V=u(G_0)$ on each uniformization neighbourhood, where $G_0\subset\D$ is the null set outside the compact Lipschitz pieces of $u$ \cite[Proposition~2.6 and Theorem~1.10]{MR}.
\end{definition}
We fix a Borel representative of $X_0$. By metric differentiability, the area formula and inner regularity, there exist countably many compact sets $K_i\subset X\setminus X_0$ covering $X\setminus X_0$ up to a $\HH^2$-null set, such that $u^{-1}|_{K_i}$ is bi-Lipschitz for each $i$. Here and below $u$ denotes the corresponding local uniformization. 
	
For each $i$, let $\varphi_i:X\to\R^2$ be a Lipschitz extension of $u^{-1}|_{K_i}$. Corollary \ref{cor:rect-charts} implies that $(K_i,\varphi_i)$ is a 2-dimensional chart whenever $\HH^2(K_i)>0$. We denote the tangent norm $|\cdot|_{i,x}=(\Phi_i^x)^*$. We set $|\cdot|_x=|\cdot|_{i,x}$ at a.e. point of the corresponding piece.  If $\mathcal H^2(X_0)>0$ then \cite[Theorem 1.8]{MR} and \cite[Proposition 4.4]{ES} give the zero-dimensional chart $(X_0,\mathbf 0)$, hence the pointwise tangent space on $X_0$  is $\R^0$, with its trivial norm.
The resulting isometry class of $|\cdot|_x$ is independent of the chosen charts \cite[Lemma 5.2]{ES}.}

We fix this setup and notation for the remainder of this section. The first lemma rules out the existence of analytically 1-dimensional parts of metric surfaces.

\begin{lemma}\label{lem:no-1-dim-charts}
If $(U,\varphi)$ is a chart of $(X,d,\HH^2)$, then either $\dim(U,\varphi)=0$ or $\dim(U,\varphi)=2$. In the first case, we have $\HH^2(U\setminus X_0)=0$; in the second case $\HH^2(X_0\cap U)=0$ and $U$ is countably 2-rectifiable.
%on a metric surface has dimension $0$ or $2$. Modulo $\HH^2$-null sets, in the former case $U\subset X_0$; in the latter, $U\subset X\setminus X_0$ is $2$-rectifiable.
\end{lemma}
\begin{proof}

Suppose $(U,\varphi)$ is a chart of $(X,d,\HH^2)$ with $\HH^2(X_0\cap U)=0$. Then $U$ can be covered up to a $\HH^2$-null set by countably many of the sets $K_i$. It follows that $U$ is 2-rectifiable and moreover $\dim(U,\varphi)=2$ (cf. \cite[Lemma 5.2]{ES}). 

If $(U,\varphi)$ is a chart of $(X,d,\HH^2)$ with $\HH^2(U\cap X_0)>0$ then -- since $(X_0,\mathbf{0})$ is a zero-dimensional chart \cite[Theorem 1.8]{MR} -- we have that $\dim(U,\varphi)=0$ \cite[Lemma 5.2]{ES}. Moreover we have $\HH^2(U\setminus X_0)=0$, since otherwise $U$ would intersect one of the 2-dimensional charts $(K_i,\varphi_i)$ on a set of positive measure.
\end{proof}

\begin{proof}[Proof of Theorem \ref{thm:main2}]
Assume that $X$ satisfies the property (ET).
Lemma \ref{lem:key-lemma} applied to  $\mu=\HH^2$ and the uniformizing map $u$ implies that $|\cdot|_{i,x}=\apmd_{\varphi_i(x)}u$ $\HH^2$-a.e. $x\in K_i$. Since $X$ satisfies (ET), it follows that the pointwise norm is Euclidean for $\HH^2$-a.e. $x\in X\setminus X_0$. On the other hand for $\HH^2$-a.e. $x\in X_0$ the pointwise norm is trivial since $X_0$ is a 0-dimensional chart. Thus $(X,d,\HH^2)$ is infinitesimally Hilbertian by \cite[Corollary 6.7(ii)]{ES}.

In the opposite direction, suppose $(X,d,\HH^2)$ is infinitesimally Hilbertian, and let $f\in W^{1,2}(\D,X)$ be a Sobolev map. For a.e. $z$ with $\J f(z)=0 $ we have that $\apmd_zf$ is a degenerate seminorm on $\R^2$, hence induced by a positive semidefinite quadratic form. The set $\{z: \J f(z)>0\}\subset \D$ admits a partition up to a null-set $\{K_i\}_i$ into compact subsets of $\D$ so that $f|_{K_i}$ is bi-Lipschitz. Thus $(f(K_i),g_i)$ are 2-dimensional charts, where $g_i$ is a Lipschitz extension of $f^{-1}|_{f(K_i)}$ (Corollary \ref{cor:rect-charts}). Denote the pointwise norm given by this chart by $|\cdot|_i$. By Lemma \ref{lem:key-lemma} and the fact that $f_\ast(\mathcal L^2|_{K_i})\ll\HH^2$ we have
\[
|\cdot|_{i,f(z)}=\apmd_{z}f\quad \mathcal L^2\textrm{-a.e.}\ z\in K_i.
\]
Since $(X,d,\HH^2)$ is infinitesimally Hilbertian the pointwise norms $|\cdot|_{i,f(z)}$ are Euclidean for a.e. $z\in K_i$. Thus $\apmd_zf$ is a Euclidean  seminorm for a.e. $z\in \D$, proving that $X$ has the property (ET).
\end{proof}

The problem with Theorem \ref{thm:main}(iii) for arbitrary metric surfaces is that we cannot relate the (norm of the) approximate tangent planes to the pointwise norms given by the charts at points in the singular part of the surface. Nevertheless, we can still relate them for points outside the singular part.

\begin{theorem}\label{thm:rect-part}
Let $X$ be a metric surface and $E\subset X$ a $\HH^2$-measurable set with  $\HH^2(E\cap X_0)=0$. Then $E$ is countably 2-rectifiable and, for $\HH^2$-a.e. $x\in E$, the approximate tangent of $E$ at $x$ is isometric to $(\R^2,|\cdot|_x)$.
\end{theorem}

\begin{proof}
The compact pieces $K_i$ cover $X\setminus X_0$ up to a null set and are bi-Lipschitz images of subsets of $\R^2$. Thus $E$ is countably 2-rectifiable.  For a.e. $x\in E\cap K_i$, put $z=\varphi_i(x)$ and $s_z=\apmd_z u$ for the corresponding local parametrization. Since $u^{-1}(E\cap K_i)$ has Lebesgue density one for a.e.  $z\in u^{-1}(E\cap K_i)$, the restriction $u|_{u^{-1}(E\cap K_i)}$ satisfies $\md_z(u|_{u^{-1}(E\cap K_i)})=s_z$ for a.e. $z\in u^{-1}(E\cap K_i)$. The approximate tangent construction for rectifiable sets \cite[Section 5]{AK} identifies the tangent of $E$ for $\HH^2$-a.e. $x\in E\cap K_i$ with $(\R^2,s_{\varphi_i(x)})$, up to a linear isometry. Since $s_{\varphi_i(x)}=|\cdot|_x$ for $\HH^2$-a.e. $x\in E\cap K_i$ by Lemma \ref{lem:key-lemma}, this proves the theorem.
\end{proof}

\subsection{Metric length surfaces}
In this subsection we assume that $X$ is a metric length surface. We continue using the setup and notation from the beginning of Section \ref{sec:metr-surf}.

\begin{proposition}\label{prop:unrect}
Let $X$ be a metric length surface. Then the singular part $X_0$ of $X$ is purely 2-unrectifiable.
\end{proposition}

\begin{proof}
We work locally in a uniformization neighbourhood $V\subset X$. Let $V_1$ denote its curve-regular points: these are the interior points of arc-length parametrized curves in $V$ at which the curve is strongly metrically differentiable with metric derivative one. By \cite[Proposition 5.3 and Theorem 1.10]{MR} we have
\[
\HH^2(X_0\cap V_1)=0.
\]
Choose a Borel $\HH^2$-null set $N\subset V$ containing $X_0\cap V_1$. Every rectifiable curve in $V$ has image contained in $V_1$ up to a $\HH^1$-null set, by strong metric differentiability of its arc-length parametrization \cite[Theorem 2.4]{MR}. Curve-regularity is local, so the same statement holds for the portion in $V$ of any rectifiable curve in $X$.

Let $\alpha:C\to X$ be bi-Lipschitz, with $C\subset\R^2$ compact, and let $L$ be a Lipschitz constant of $\alpha$. For each horizontal line, write $C_t=\{s:(s,t)\in C\}$. If $C_t$ has at least two points, the fragment $s\mapsto\alpha(s,t)$ extends to a $2L$-Lipschitz curve $\gamma_t$ defined on the smallest interval $I_t$ containing $C_t$. Indeed, on each interval $(a,b)\subset I_t\setminus C_t$, the length-space assumption allows us to join $\alpha(a,t)$ to $\alpha(b,t)$ by a curve of length at most $2L(b-a)$; parametrizing with constant speed gives the extension. Clearly $\gamma_t$ is a $2L$-Lipschitz curve. If $C_t$ has at most one point there is no contribution to one-dimensional measure.

The image of $\gamma_t$ meets $(X_0\cap V)\setminus N$ in a $\HH^1$-null set, since $(X_0\cap V)\setminus N$ contains no curve-regular point of $V$. Since $\alpha$ is bi-Lipschitz we therefore have
\[
\LL^1\{s\in C_t:\alpha(s,t)\in(X_0\cap V)\setminus N\}=0
\]
for every $t$. Fubini's theorem yields
\[
\LL^2\bigl(C\cap\alpha^{-1}((X_0\cap V)\setminus N)\bigr)=0.
\]
Thus $\HH^2((\alpha(C)\cap X_0\cap V)\setminus N)=0$ and, since $\HH^2(N)=0$, we have $\HH^2(\alpha(C)\cap X_0\cap V)=0$. Covering $X$ by countably many uniformization neighbourhoods $V$, this yields $\HH^2(X_0\cap \alpha(C))=0$. Since 2-rectifiable sets can be decomposed up to a $\HH^2$-null set into bi-Lipschitz images of planar sets \cite{K94}, this implies that $X_0$ is purely 2-unrectifiable.
\end{proof}

\begin{proof}[Proof of Theorem \ref{thm:main}]
	The equivalence of (i) and (ii) is proved in Theorem \ref{thm:main2}. The implication (iii)$\implies$(i) is  true in all metric spaces, see Section \ref{sec:general-relations}.  It remains to verify  the implication  (ii)$\implies$(iii). 
	
	Let $E\subset X$ be 2-rectifiable.  We write $\operatorname{Tan}(E,x)$ for its approximate tangent plane at $x$. Then $\HH^2(E\cap X_0)=0$ by Proposition \ref{prop:unrect}. Consequently Theorem \ref{thm:rect-part} implies that 
	\begin{align}\label{eq:tan}
		\operatorname{Tan}(E,x) \cong (\R^2,|\cdot|_x)\quad\textrm{for}\quad\HH^2\textrm{-a.e.}\ x\in E.
	\end{align}
	If $(X,d,\HH^2)$ is infinitesimally Hilbertian, then $\operatorname{Tan}(E,x)$ is Euclidean for $\HH^2$-a.e. $x\in E$, due to  \eqref{eq:tan}.  This verifies  (ii)$\implies$(iii) and
 completes the proof of Theorem \ref{thm:main}.
\end{proof}

\subsection{Universal infinitesimal Hilbertianity}\label{sec:universal}
In this subsection we prove Theorem \ref{thm:main3}. When several measures are under consideration, we write $|D_\eta f|$, $\Mod_{2,\eta}$,  and $|\cdot|_{\eta,x}$ for the minimal weak upper gradient, 2-modulus, and pointwise norm, taken with respect to the measure $\eta$.

\begin{lemma}\label{lem:no-2dim-charts-on-sing-part}
Let $X$ be a metric surface and $\mu$ a locally finite measure on $X$. If $(U,\varphi)$ is a 2-dimensional chart of $(X,d,\mu)$, then $\mu|_U\ll\HH^2|_U$ and  $\mu(U\cap X_0)=0$.
\end{lemma}
\begin{proof}
By the proof of \cite[Proposition 4.13]{ES} we have that $\varphi_\ast(\mu|_U)\ll \LL^2$. Since $\varphi$ maps $\HH^2$-null sets to $\LL^2$-null sets, this implies that $\mu|_U\ll \HH^2|_U$.

Let  $X_1$ denote the curve regular points of $X$; see \cite[Definition 5.1]{MR}. By \cite[Proposition 5.3]{MR} there exists a Borel set $N\subset X_0$ containing $X_0\cap X_1$ with $\HH^2(N)=0$. For any rectifiable curve $\gamma$ parametrized by arclength we have $|\gamma^{-1}(X_0\setminus N)|=0$, since $\gamma(t)\in X_1$ for a.e. $t$. It follows that 
\begin{align}\label{eq:mod-zero}
\Mod_{2,\mu}\Gamma^+_{X_0\setminus N}=0.
\end{align}
where $\Gamma^+_A=\big\{\gamma:\ \int_\gamma\chi_A\ud s>0 \big\}$. On the other hand $\mu(U\cap N)=\HH^2(U\cap N)=0$ by the absolute continuity $\mu|_U\ll\HH^2|_U$, implying $\Mod_{2,\mu}\Gamma^+_{U\cap N}=0$. Together with \eqref{eq:mod-zero} this yields $\Mod_{2,\mu}\Gamma^+_{X_0\cap U}=0$. If $\mu(X_0\cap U)>0$, $(X_0\cap U,\mathbf 0)$ would be a 0-dimensional chart \cite[Proposition 4.4]{ES}, contradicting the 2-dimensionality of $(U\cap X_0,\varphi)$. Thus $\mu(X_0\cap U)=0$.
\end{proof}

\begin{proof}[Proof of Theorem \ref{thm:main3}]
By Theorem \ref{thm:main2}, we may assume that $X$ has property (ET). Define $\nu:=\mu+\HH^2$, and let $(U,\varphi)$ be a chart of $(X,d,\nu)$. If $\dim(U,\varphi)\le 1$, the pointwise norm is automatically Hilbertian; assume therefore that $\dim(U,\varphi)=2$. By Lemma \ref{lem:no-2dim-charts-on-sing-part} we have $\nu(X_0\cap U)=0$ and $\nu|_U\ll \HH^2|_U$. This and Lemma \ref{lem:key-lemma} applied to $\nu$ give that
\begin{align*}
	(\R^2,|\cdot|_{\nu,x})\stackrel{\textrm{isometr.}}{\simeq}(\R^2,\apmd_{\varphi_i(x)}u_i)\quad \nu-\textrm{a.e.}\ x\in U,
\end{align*}
where $u_i$ is the corresponding local uniformization map and $\varphi_i$ a Lipschitz extension of $u_i^{-1}|_{K_i}$. Since $X$ has property (ET), it follows that $|\cdot|_{\nu,x}$ is Hilbertian $\nu$-a.e. $x\in U$. This shows that $(X,d,\nu)$ is infinitesimally Hilbertian. 
	
We claim that $(X,d,\mu)$ is infinitesimally Hilbertian. By suitable localization arguments we may assume that the measures $\mu$ and $\nu$ are finite. Since $\mu\le \nu$, the minimal weak upper gradients with respect to $\mu$ and $\nu$, respectively, satisfy $|Df|_\mu\le |Df|_\nu$ $\mu$-almost everywhere. Indeed, if $f\in L^2_\loc(\nu)$ has $\nu$-weak upper gradient $g\in L^2_\loc(\nu)$, then -- since $\Mod_{2,\mu}\le \Mod_{2,\nu}$ -- the upper gradient inequality $\displaystyle |f(\gamma(b))-f(\gamma(0))|\le \int_\gamma g\ud s$ holds for $\Mod_{2,\mu}$-a.e. curve, implying $|Df|_\mu\le g$. In particular, $|Df|_\nu$ is well-defined up to a $\mu$-null set. Note that, for Lipschitz $f$ we have
	\begin{align}\label{eq:energy-ineq}
		\int|Df|_\mu^2\ud\mu\le \int |Df|_\nu^2\ud\mu\le \int(\operatorname{Lip} f)^2\ud\mu.
	\end{align}
Here $\operatorname{Lip}f(x)=\limsup_{y\to x}\frac{|f(y)-f(x)|}{d(y,x)}$ is the pointwise Lipschitz constant of $f$. Denote $Q(f):=\int|Df|_\nu^2\ud\mu$ for bounded Lipschitz functions $f\in L^2(\mu)$, and observe that $Q$ is quadratic since $(X,d,\nu)$ is infinitesimally Hilbertian. By energy density of Lipschitz functions \cite{AGS13,EB23} and lower semicontinuity of $f\mapsto \int|Df|_\mu^2\ud\mu$ with respect to $L^2(\mu)$-convergence, \eqref{eq:energy-ineq} yields
	\begin{align*}
		\int|Df|_\mu^2\ud\mu=\overline Q(f), \quad f\in W^{1,2}(X,d,\mu),
	\end{align*}
	where $\overline Q$ denotes the $L^2(\mu)$-relaxation of $Q$. Since the relaxation of a quadratic functional is quadratic, it follows that $f\mapsto \int|Df|_\mu^2\ud\mu$ is quadratic and thus $W^{1,2}(X,d,\mu)$ is a Hilbert space. This completes the proof.
\end{proof}

\begin{remark}
The preceding proof actually yields the following statement which is of independent interest: Let $(X,d)$ be a locally complete separable metric space and let $\mu\le \nu$ be two locally finite Borel regular measures on $X$. If $(X,d,\nu)$ is infinitesimally Hilbertian, then $(X,d,\mu)$ is infinitesimally Hilbertian.
\end{remark}

\end{document}